\documentclass{svjour3} 
\smartqed  

\usepackage{amsfonts}
\usepackage{amsmath} 
\usepackage{amssymb}  

\usepackage{graphics}
\usepackage{multirow}
\usepackage{epsfig}
\usepackage[ruled,vlined]{algorithm2e}
\usepackage{subfig}
\usepackage{amsfonts}
\usepackage{mathrsfs}
\usepackage{bbm}
\usepackage{color}
\usepackage{pgfpages}
\usepackage{tikz}
\usepackage{tkz-berge}
\usetikzlibrary{arrows,snakes,backgrounds}                    

\smartqed  
\usepackage{graphicx}

\usepackage[english]{babel}                 	
\usepackage[T1]{fontenc}                    	
\usepackage{ae}                             	
\usepackage{setspace}      										
\usepackage{longtable}  
\usepackage{algpseudocode}   										
\usepackage{graphicx}      										
\usepackage{tikz}          										
\usepackage[permil]{overpic}       						
\usepackage{fancybox}													
\usepackage{ifthen}														
\usepackage{tikz}
\usepackage{tikz-3dplot}
\usepackage[numbers]{natbib}
\usepackage{fancyhdr, graphicx}
\usepackage{pgfgantt}
\usepackage{bm}
\SetKwComment{Comment}{/* }{ */}
\renewcommand\familydefault\sfdefault
\let\xell\ell
\usepackage[symbolgreek,defaultmathsizes]{mathastext}
\let\ell\xell
\renewcommand\familydefault\rmdefault
\DeclareMathOperator*{\argmin}{arg\,min}
\newcommand{\keff}{k_{\text{eff}}}

\newcommand{\gphi}{\textbf{\textbf{$\phi_g$}}}
\newcommand{\grho}{\textbf{\textbf{$\rho$}}}

\begin{document}

\title{A multi-fidelity adaptive dynamical low-rank based optimization algorithm for fission criticality problems
}
\subtitle{}


\author{Carmela Scalone         \and
        Lukas Einkemmer  \and Jonas Kusch \and Ryan G. McClarren
}


\institute{C. Scalone \at
              University of L'Aquila, Italy \\
            \email{carmela.scalone@univaq.it} 
                       \and
          L. Einkemmer \at
              University of Innsbruck, Austria \\
             \email{lukas.einkemmer@uibk.ac.at} 
             \and
                           J. Kusch \at
            Norwegian University of Life Science, Norway \\
                            \email{jonas.kusch@nmbu.no} 
                            \and
             R. J. McClarren \at
              University of Notre Dame, USA\\
              \email{rmcclarr@nd.edu}}

\date{Received: date / Accepted: date}

\maketitle

\begin{abstract}
Computing the dominant eigenvalue is important in nuclear systems as it determines the stability of the system (i.e. whether the system is sub or supercritical). Recently, the work of Kusch, Whewell, McClarren and Frank \cite{KWMF} showed that performing a low-rank approximation can be very effective in reducing the high memory requirement and computational cost of such problems. In this work, we propose a rank adaptive approach that changes the rank during the inverse power iteration. This allows us to progressively increase the rank (i.e. changing the fidelity of the model) as we get closer to convergence, thereby further reducing computational cost. We then exploit this multi-fidelity approach to optimize a simplified nuclear reactor. In this case the system is parameterized and the values of the parameters that give criticality are sought. 
\end{abstract}

\section*{Introduction}
\label{intro}
Recently energy production in fission reactors has seen increased interest. Nuclear fission plays an important role in the energy supply chain. Hence ensuring safe conditions and optimal performance of nuclear systems has significant  environmental and economic relevance. An important aim of nuclear engineering is to achieve a sustainable chain reaction where the neutrons produced by fission balance the neutrons that are either absorbed or leave the reactor, thus ensuring safety and optimal performance. Neutron transport theory  describes
the process to maintain a stable and continuous fission reaction, so accurate mathematical modelling and obtaining efficient solutions of the resulting problems is important.

In this paper, we are interested in the study of the criticality of a (simplified) reactor, see \cite{McClarrenBook, M,StaceyBook}. 
Mathematically this is formulated as a generalised eigenvalue problem. In particular, we focus on the numerical computation of the dominant eigenvalue, called the \textit{effective eigenvalue} (denoted by $\keff$), whose magnitude describes  the distance from \textit{criticality}, i.e. from the equilibrium state of the reactor where the loss and the gain of neutrons are precisely balanced. 
In particular, the system is in the desired \textit{critical} state when $\keff = 1$, while for $\keff>1$ (corresponding to a runaway chain reaction) and $\keff<1$  (corresponding to an exponential decay of the neutrons available for fission) we have \textit{supercritical} and \textit{subcritical} states, respectively.
An efficient and accurate numerical solution of the effective eigenvalue is important, since it can be used to design and optimize nuclear reactors. Often $\keff$ needs to be determined to a stringent accuracy, as even small deviations from $1$ can result in a significant growth or decay of the neutron population over time.

The mathematical description for neutron transport used in this work is a multigroup diffusion model. In this model an elliptic diffusion operator is used to model the transport of neutrons and the collision of neutrons with each other and the material of the reactor. The neutron energies are packed into finite energy ranges known as energy groups. This model can be obtained from kinetic theory in the limit of (relatively) strong collisionality. In order to solve this PDE a number of techniques has been proposed in the literature. Once a space discretization has been performed, the resulting generalized eigenvalue problem for $\keff$ can be treated by techniques from linear algebra such as (inverse) power or Krylov iteration, see \cite{KPN,Ortega,WPK}  and reference therein.

Since such simulations can be computationally expensive, especially in a multi-query context (i.e. for optimization or uncertainty quantification), model order reduction has been investigated in the literature, see \cite{Buchan,GR,PM,PMF,Sun}. Dynamical low rank approximation (DLRA) is an approach to find an approximation of the solution of a time-dependent matrix differential equation on a low rank manifold, see \cite{KL}. Several suitable integrators (see,e .g., \cite{LO,CL,CKL, HNS}) to provide a numerical solution have been developed and  DLRA has been demonstrated to be efficient for a range of problems including quantum physics (see, e.g., \cite{lrqm1,lrqm2,lrqm3,lrqm4}), kinetic equations (see, e.g., \cite{lrkin1,lrkin2, PMF,lrkin4,lrkin5,lrap4}), biology \cite{lrbiology1,lrbiology2,lrbiology3}, and uncertainty quantification \cite{uncertainty1}. Much progress has also been made with regards to the development of structure preserving algorithms (see, e.g. \cite{lrcons1,EOS,PM,lrcons4}), asymptotic preserving schemes (see, e.g., \cite{lrap1,lrap2,lrap3}) and rank adaptive algorithms (see, e.g., \cite{CKL, HNS,lradap3,HS}).
In \cite{GKS}, a DRLA approach for computing rightmost eigenpairs of linear operators has been provided. Recently,  it has been shown that DLRA can be also used in the context of an inverse power iteration, see \cite{KWMF}. Essentially, each iteration of the method corresponds to a time step in the classic DLRA algorithm.

An additional benefit of this approach, besides the overall reduction in memory and computational cost, is that the resulting equations for the DLRA inverse power method require only linear solves of size $N_x$ (the number of points in the spatial discretization) and $G$ (the number of energy groups), but not the solution of a linear system of size $N_x \cdot G \times N_x \cdot G$. This makes it feasible to use direct solvers for these smaller linear systems and still being able to obtain an efficient numerical method.

In this paper, we propose a specific rank adaptive algorithm for computing the effective eigenvalue using a dynamical low-rank power iteration. The aim is to obtain an algorithm capable of providing an accurate approximation of the eigenvalue and the eigenmatrix (i.e. keeping the low-rank truncation error as small as possible), while reducing the required computational cost. To accomplish this we perform the intermediate iterations with lower and gradually increasing rank until a certain accuracy is achieved. The stationarity of the power iteration is used as an indicator for increasing the rank.
Rank adaptation techniques have proven to be particularly efficient in several DLRA frameworks, see \cite{CKL,EOS, HS, HNS}.
In practice it is often useful to optimize the parameters of a design in order to get a $\keff$ that is as close as possible to criticality. In this paper, we therefore also propose a low rank based optimization algorithm. We illustrate the effectiveness of this algorithm by performing numerical experiments for a sphere with a hollow center, filled by uranium and with two steel casings of different thickness. The values of the thickness are the parameters that are optimized. We also combine the optimization algorithm with the rank adaptive algorithm developed here.

The paper is organized as follows: in Section 1 we recap the mathematical model used, the power method and the dynamical low rank power method of \cite{KWMF}, which is the starting point of the present paper. In Section 2 we present the new rank adaptation technique and related numerical experiments are provided in Section 3. Section 4 is dedicated to the low rank optimization techniques and we conclude with numerical experiments in Section 5. 
\section{Framework}
In this section, we present the problem and the notation used. In particular, following \cite{KWMF} we introduce the mathematical model and the related compact matrix notation, obtained after a numerical discretization in space has been performed.

The multigroup diffusion equation is given by 
\begin{equation}
\label{problem}
-\nabla \cdot (D_g(r) \nabla \phi_g(r)) + \Sigma_{t,g}(r) \phi_g(r) =  \frac{\chi_g}{k_{\text{eff}}} \sum_{g'} \nu \Sigma_{f,g'}(r) \phi_g'(r) + \sum_{g'} \Sigma_{s,g',g}(r) \phi_{g^{'}}(r)
\end{equation}
where $r \in \Omega \in \mathbb{R}^d$ is the spatial domain, $\Sigma_{t,g}(r)$  the total cross section of the energy group $g$ at a spatial position $r$, $\Sigma_{f,g}(r)$  the fission cross-section and $\Sigma_{s,g', g}(r)$ is the scattering cross section between groups $g$ and $g'$, $\nu$  is the mean number of particles produced per fission event and $\chi_g$  is the fission scalar flux over the energy range of group $g$ at position $r$ corresponding to the maximal eigenvalue $\keff$.
The aim is to compute $\phi_g(r)$, which is the integral of the scalar flux over the energy range of group $g$ at position $r$, corresponding to the maximum eigenvalue $\keff$.
To introduce a numerical discretization, we consider a spatial grid $r_1, r_2, ..., r_{N_x}$ and energy groups $g \in \lbrace1, ..., G\rbrace$.
Therefore,  we represent the numerical solution as a matrix  $\phi = \lbrace\phi _{jg}\rbrace_{j,g=1}^{N_x,G}$, where $\phi_{jg} = \phi_g(r_j)$.
We now review the steps to arrive at a formulation of the problem as a matrix equation.
The diffusion coefficient can be written as
\begin{align}\label{eq:expansionD}
    D_g(r) = \sum_{\ell = 1}^{N_m} \rho_{\ell}(r)D_{g}^{(\ell)},
\end{align}
where $N_m$ is the number of materials and $D_{g}^{(\ell)}$ is the diffusion coefficient (or any other material coefficient) for material $\ell$ and energy group $g$. The functions $\rho_{\ell}(r)$ denotes the density of material $\ell$. 
The evaluation of $D_g$ at the cell interfaces $r_{j+1/2} = (r_j+r_{j+1})/2$, need to approximate the material coefficient though its harmonic mean, i.e.,
\begin{align*}
    D_g(r_{j+1/2}) =: D_{g,j+1/2} = 2\frac{D_{g,j} D_{g,j+1}}{D_{g,j} + D_{g,j+1}} = 2\frac{\sum_{\ell,k}\rho_j^{(\ell)}D_g^{(\ell)}\rho_{j+1}^{(k)}D_g^{(k)}}{\sum_{\ell}(\rho_j^{(\ell)}+\rho_{j+1}^{(\ell)})D_g^{(\ell)}}.
\end{align*}
Noting that only two terms in the sums can be non-zero, we get
\begin{align}\label{eq:approxD}
    D_{g,j+1/2} = \sum_{\ell,k=1}^{N_m}\rho_j^{(\ell)}\rho_{j+1}^{(k)}\left(\rho_j^{(\ell)}+\rho_{j+1}^{(k)}\right)\frac{D_g^{(\ell)}D_g^{(k)}}{D_g^{(\ell)}+D_g^{(k)}}.
\end{align}
The diffusion operator in one dimension is discretized as follows
\begin{align*}
    \nabla \cdot D_g(r)\nabla\phi_g(r)\Big|_{r=r_j} \approx \left( \mathbf{D}(g) \gphi\right)_j
\end{align*}
where $ \gphi \in\mathbb{R}^{N_x}$ collects the scalar flux at all spatial cells. The matrix $\mathbf{D}(g)\in\mathbb{R}^{N_x\times N_x}$ has values
\begin{align*}
    D_{j,j\pm 1}(g) &= \pm\frac{1}{\Delta r\cdot V_j} D_{g,j\pm 1/2}S_{j\pm 1/2}, \\
    D_{j,j}(g) &= -\frac{1}{\Delta r\cdot V_j} \left[D_{g,j+1/2}S_{j+1/2}+D_{g,j-1/2}S_{j-1/2}\right],
    \end{align*}
where $\Delta r$ is the size of each radial element. The surface area between cell $j$ and $j\pm 1$ is denoted by $S_{j\pm 1/2}$ and the area of cell $j$ is denoted by $V_j$. The choice of these terms defines the spatial geometry. In our numerical experiments, we look at spherical domains, where we have $V_j = \frac{4\pi}{3}(r_{i+1/2}^3-r_{i-1/2}^3)$ and $S_{j\pm 1/2}=4\pi r_{i\pm 1/2}^2$.
Using \eqref{eq:approxD} in the above definition of $ \mathbf{D}(g)$, lets us write ${\mathbf{D}}(g){\phi}_g$ as

\begin{align*}
  {\mathbf{D}}(g){\gphi} = \sum_{\ell,k = 1}^{N_m} \frac{D_g^{(\ell)}D_g^{(k)}}{D_g^{(\ell)}+D_g^{(k)}} \mathbf{D}^{(\ell,k)} \gphi,
\end{align*}
where we use
\begin{align*}
    D^{(\ell,k)}_{j,j\pm 1} &= \pm\frac{\rho_j^{(\ell)}\rho_{j\pm 1}^{(k)}}{\Delta r\cdot V_j}  (\rho_{\ell}(r_{j})+\rho_{k}(r_{j\pm 1}))S_{j\pm 1/2}, \\
    D^{(\ell,k)}_{j,j} &= -\frac{1}{\Delta r\cdot V_j} \left[ \rho_j^{(\ell)}\rho_{j+1}^{(k)}\left(\rho_j^{(\ell)}+\rho_{j+1}^{(k)}\right)S_{j+1/2}+ \rho_j^{(\ell)}\rho_{j-1}^{(k)}\left(\rho_j^{(\ell)}+\rho_{j-1}^{(k)}\right)S_{j-1/2}\right].
\end{align*}

The diagonal matrix ${\grho}^{(\ell)}\in\mathbb{R}^{N_x\times N_x}$ with entries $\rho^{(\ell)}_{j j} = \rho_{\ell}(r_j)$ can be used to write
\begin{align*}
   \frac{\chi_g \nu}{k}\sum_{g'}\Sigma_{f,g'}(r)\phi_{g'}(r) &= \frac{\chi_g\nu}{k}\sum_{g',\ell} \rho_{\ell}(r_j)\Sigma_{f,g'}^{(\ell)}\phi_{g'}(r_j) = \frac{\chi_g\nu}{k}\sum_{g',\ell} \Sigma_{f,g'}^{(\ell)}{\grho}^{(\ell)} \mathbf{\phi}_{g'}, \\
   \sum_{g'}\Sigma_{s,g',g}(r_j) \phi_g(r_j) &= \sum_{g',\ell}\rho_{\ell}(r_j)\Sigma_{s,g',g}^{(\ell)} \phi_g(r)= \sum_{g',\ell}\Sigma_{s,g',g}^{(\ell)}{\mathbf{\rho}}^{(\ell)} \mathbf{\phi}_{g'}.
\end{align*}

For a given group $g$, the diffusion equation reads
\begin{align*}
  -\sum_{\ell,k = 1}^{N_m} \frac{D_g^{(\ell)}D_g^{(k)}}{D_g^{(\ell)}+D_g^{(k)}} D^{(\ell,k)} \phi_g + \sum_{\ell = 1}^{N_m}\Sigma_{t,g}^{(\ell)}{\rho}^{(\ell)}\phi_g = \frac{\chi_g\nu}{\keff}\sum_{g',\ell} \Sigma_{f,g'}^{(\ell)}{\mathbf{\rho}}^{(\ell)} \mathbf{\phi}_{g'}+\sum_{g',\ell}\Sigma_{s,g',g}^{(\ell)}{\mathbf{\rho}}^{(\ell)} \mathbf{\phi}_{g'}.
\end{align*}
We define ${\widetilde{\mathbf{\Sigma}}}_{f}^{(\ell)}= \left(\chi_g\nu\Sigma_{f,g'}^{(\ell)}\right)_{g,g'=1}^G$ and the diagonal matrices $ \mathbf{M}^{(\ell,k)}\in\mathbb{R}^{G\times G}$ with entries $M_{gg}^{(\ell,k)} = D_g^{(\ell,k)}$ as well as $\mathbf{\Sigma}_t^{(\ell)}\in\mathbb{R}^{G\times G}$ with $\Sigma_{t,gg}^{(\ell)}=\Sigma_{t,g}^{(\ell)}$. Then we have the following matrix equation for $\mathbf{\phi} \in \mathbb{R}^{N_x \times G}$
\begin{align*}
    -\sum_{\ell,k}  \mathbf{D}^{(\ell,k)} \mathbf{\phi}  \mathbf{M}^{(\ell,k)} + \sum_{\ell}     \mathbf{{\rho}}^{(\ell)}\mathbf{\phi}\mathbf{{\Sigma}}_t^{(\ell)} = \frac{1}{\keff}\sum_{\ell} \mathbf{\rho}^{(\ell)} \mathbf{ \phi}{\widetilde\Sigma}_{f}^{(\ell)}+\sum_{\ell}{\rho}^{(\ell)}\mathbf{\phi} \mathbf{\Sigma}_{s}^{(\ell)}.
\end{align*}
Assigning an iteration index to $\phi$ and setting $\widetilde{\phi}^{n+1} = \keff \phi^{n+1}$, we can write
\begin{equation}
\label{eq_powerit}
- \sum_{l,k}   \mathbf{D}^{(\ell,k)} \widetilde{\mathbf{\phi}}^{n+1}  \mathbf{M}^{(\ell,k)}+ \sum_l   \mathbf{\rho}^{(\ell)} \widetilde{\mathbf{\phi}}^{n+1}
\left(  \mathbf{\Sigma}_t^{(\ell)} -  \mathbf{\Sigma}_s^{(\ell)}\right) = \sum_l  \mathbf{\rho}^{(\ell)} {\mathbf{\phi}}^n   \widetilde{\mathbf{\Sigma}}_f^{(\ell)}.
\end{equation}
Note that this is an implicit update. That is, in order to obtain $\tilde{\mathbf{\phi}}^{n+1}$ from $\mathbf{\phi}^n$ we need to solve a linear system of size $N_x \times G$.  The power method to compute the dominant eigenpair is then given as follows 
\begin{enumerate}
\item Start with a normalized initial guess $\phi_g^0$ with $g = 1,\ldots ,G$ and set $n=0$. 
\item Compute $\widetilde{\phi}^{n+1}$ from \eqref{eq_powerit} using $\phi^n$.
\item Set $k_{n+1} = \| \widetilde{\phi}^{n+1} \|$ and $\phi^{n+1} = \widetilde{\phi}^{n+1} /k_{n+1}$.
\item If $|k_{n+1}-k_n| < tol$ stop
\item Else set $n = n+1$ and repeat from Step 2.
\end{enumerate}

\subsection{A low rank power method}

The key ideas consist of interpreting the $n$-th iteration of the power method as a timestep of the unconventional integrator of \cite{CL}.
We denote by $k_n$ the approximation of $\keff$ at the $n$-th iteration  and  consider a rank $r$ representation of $\phi^n$, i.e. $$\phi^n= X^n S^n W^{n,T},$$ where $X^n \in \mathbb{R}^{n_x \times r}, S^n \in \mathbb{R}^{r \times r}$ and $W^n \in \mathbb{R}^{G \times r}$.
According to the dynamical low rank approximation \cite{KL}, $X^n \in \mathbb{R}^{N_x \times r}$ and $W^n \in \mathbb{R}^{G \times r}$ have orthonormal columns and $ S^n \in \mathbb{R}^{r \times r}$ is non singular. Increasing the rank $r$ increases accuracy but also computational cost.
Defining $\Sigma_{\ell} = \Sigma^{(\ell)}_t-\Sigma_s^{(\ell)}$ and given an accuracy parameter $\theta$, the dynamical low rank based power method of \cite{KWMF}, for fixed rank $r$, takes the form
\begin{enumerate}
\item \textbf{K-step}: We define $K^n = X^n S^n$. Update $X^n$ to $X^{n+1}$ via
$$- \sum_{\ell,k}D^{(\ell,k)} K^{n+1} \widehat{M}_n^{(\ell,k)}+ \sum_{\ell} \rho^{(\ell)} K^{n+1}\widehat{ \Sigma}_n^{(\ell)} =  \sum_{\ell} \rho^{(\ell)} K^{n}\widehat{\Sigma}_{f,n}^{(\ell)},$$
where the terms $\hat{\Sigma}^{(\ell)}_{f,n} = W^{n,T} {\Sigma}_{(\ell)}W^n$, $\hat{M}^{(\ell,k)} = W^{n,T} M^{(\ell,k)} W^n$, $\hat{\Sigma}^{(\ell)}_{n} = W^{n,T} \tilde{\Sigma}^{(\ell)}_{f}W^n$, are computed in $\mathcal{O}(r^2 G^2)$ operations.
Determine $X^{n+1}$ with $K^{n+1} = X^{n+1}R$ and store $M_x = X^{n+1, T}X^n$.

\item \textbf{L-step}: We define $L^n = W^n S^{n,T}$. Update $W^n$ to $W^{n+1}$ via
$$- \sum_{\ell,k}\widehat{D}_n^{(\ell,k)} L^{n+1} M^{(\ell,k)}+ \sum_{\ell} \widehat{\rho}^{(\ell)} L^{n+1} { \Sigma}^{(\ell)} =  \sum_{\ell} \widehat{\rho}^{(\ell)} L^{n}\widetilde{\Sigma}_{f}^{(\ell)}$$

where $\widehat{\rho}_n^{(\ell)} = X^{n,T} \rho^{(\ell)}X^n$, $\widehat{D}_n^{(\ell, k)} = X^{n,T}D^{(\ell,k)}X^n$ are computed by $\mathcal{O}(r^2 N_x^2)$ operations.

Determine $W^{n+1}$ with $L^{n+1} = W^{n+1}\tilde{R}$ and store $N_w = W^{n+1, T}W^n$.

\item \textbf{S-step}: Update $S^n$ to $S^{n+1}$ via
$$- \sum_{l,k}\widehat{D}_{n+1}^{(\ell,k)} \widetilde{S}^{n+1} \widehat{M}_{n+1}^{(\ell,k)}+ \sum_l \widehat{\rho}_{n+1}^{(\ell)} \widetilde{S}^{n+1}\widehat{ \Sigma}_{n+1}^{(\ell)} =  \sum_l \widehat{\rho}_{n+1}^{(\ell)} S\widehat{\Sigma}_{f, n+1}^{(\ell)}.$$
where $S= M_x S^n N_w^{T}$. 
\item Set $k_{n+1} = \| \widetilde{S}^{n+1}  \|$ and $S^{n+1} =  \widetilde{S}^{n+1} / k^{n+1}$

\item If $| k_{n+1}-k_n | \leq \theta $ stop, else set $n=n+1$ and repeat.
\end{enumerate}
At the end, we get an approximation of the desired eigenpair $\left( k_{\text{eff}},\phi \right) $, with $\phi$ of rank $r$.

\section{Rank adaptive low-rank power method}

    \begin{algorithm}
\caption{\textbf{Rank adaptive algorithm}}\label{algorithm1}
\KwData{$\phi_0 = X_0 S_0 W_0^T$, $r_0$, $\theta_0$, $\theta$, $\rho$, $\kappa$.}
\While{$\Delta_n > \theta_n $}{
$\left( k_n,\phi_r\right) _{r=r_n} = \text{\textbf{DLRP}}(\phi_0 = U_0 S_0 V_0^T , k_{n-1}, r_n, \theta_n)$\;
$\Delta_n = | k_{n} -k_{n-1}|$\; 
$r_{n+1} = r_n+ \kappa$ \Comment*[r]{Increase the rank}\
$\phi_0 = \text{\textbf{CR}}(\phi_{r_n})$ \Comment*[r]{Merge the computed solution in $\mathcal{M}_{r_{n+1}}$}
  \If{$\theta_n > \theta $}{
    $\theta_{n+1} =  \rho \theta_n $ \Comment*[r]{Reduce the tolerance}}
  $n = n+1$\;
  }
\end{algorithm}

In this section, we provide a suitable rank adaptation for the low rank power method.
The scope is to get a final accurate approximation of the maximum eigenpair $(\keff, \phi)$, exploiting as much as possible, low rank intermediate approximations.
In order to obtain a good accuracy in the final approximation, one may have to make (quite expensive) iterations with moderate to high rank. The goal is to keep the number of these iterations as low as possible.

The adaptation strategy consists of a sequence of applications of the low rank power method, each with a fixed rank, increasing the rank as more accuracy is required.
We start by considering a (small) initial rank $r_0$ and a randomly chosen initial condition $\phi_0 = X_0 S_0 W_0^T$ of rank $r_0$.

We introduce two tolerance parameters $\theta$ and $\theta_0$.
In particular, $\theta$ is the tolerance parameter to stop the entire procedure, i.e. the accuracy required for our final approximation.
The algorithm starts computing the eigenpair $(\phi_r, \keff^r)_{r=r_0}$, via the low rank power method with fixed rank $r_0$ and  accuracy $\theta_0$  in the eigenvalue.
Then, we increase the rank by an increment $\kappa$,  setting $r_1= r_0+\kappa$.
Simultaneously, we choose a a factor  $\rho<1$ to reduce the tolerance  as $\theta_1 = \rho \theta_0 $. 
We compute the next approximation $(\phi_r, \keff^r)_{r=r_1}$, via the low rank power method with rank $r_1$ and stopping tolerance $\theta_1$.
The algorithm continues in this way, until a solution $k_n \approx \keff$ is computed with the final required accuracy $\theta$.
So, at the step $n$ of the algorithm we compute $k_n$, an approximation of the desired eigenvalue with accuracy $\theta_n = \rho^n\theta_0$, and the related eigenmatrix of rank $r_n$.
Clearly, to perform the $(n+1)$-st iteration with rank $r_{n+1}$, we need an initial condition of rank $r_{n+1}$ to initialize the DLRP (dynamical low rank power method) for the next iteration.
A good starting value can be obtained by merging the solution $\phi_n$ of rank $r_n$ obtained by DLRP into the larger manifold of rank $r_{n+1}$.
We do this performing the following algorithm (CR in the Algorithm \ref{algorithm1}).
The name CR is due to fact that we perform a change of rank.
Starting from $\phi_n = X^n S^n W^{nT}$ of rank $r_n$, we do

\begin{enumerate}
\item \textbf{K-step}: Update $X^n \in \mathbb{R}^{N_x \times r_n}$ to $X^{n+1}\in \mathbb{R}^{N_x \times 2r_n}$ via
$$- \sum_{\ell,k}D^{(\ell,k)} K^{n+1} \widehat{M}_n^{(\ell,k)}+ \sum_{\ell} \rho^{(\ell)} K^{n+1}\widehat{ \Sigma}_n^{(\ell)} =  \sum_{\ell} \rho^{(\ell)} K^{n}\widehat{\Sigma}_{f,n}^{(\ell)}.$$
Determine $\widetilde{X}^{n+1}$ with $[K^{n+1}, X^n] = \widetilde{X}^{n+1}R$ and store $M = \widetilde{X}^{n+1, T}X^n \in \mathbb{R}^{2r_n \times r_n}$.
\item \textbf{L-step}: Update $W^n \in \mathbb{R}^{G \times r_n}$ to $W^{n+1}\in \mathbb{R}^{G \times 2 r_n}$ via
$$- \sum_{\ell,k}\widehat{D}_n^{(\ell,k)} L^{n+1} M^{(\ell,k)}+ \sum_l \widehat{\rho}_n^{(\ell)} L^{n+1}\widehat{ \Sigma}^{(\ell)} =  \sum_l \widehat{\rho}^{(\ell)} L^{n}\widetilde{\Sigma}_{f}^{(\ell)}.$$
Determine $\widetilde{W}^{n+1}$ with $[L^{n+1}, W^n] = \widetilde{W}^{n+1}\widetilde{R}$ and store $N = \widetilde{W}^{n+1, T}W^n \in \mathbb{R}^{2r_n \times r_n}$.
\item \textbf{S-step}: Update $S^n \in \mathbb{R}^{r_n \times r_n}$ to $\widehat{S}^{n+1}\in \mathbb{R}^{2r_n \times 2r_n}$ via
$$- \sum_{\ell,k}\widehat{D}_{n+1}^{(\ell,k)} \widehat{S}^{n+1} \widehat{M}_{n+1}^{(\ell,k)}+ \sum_{\ell} \widehat{\rho}_{n+1}^{(\ell)} \widehat{S}^{n+1}\widehat{ \Sigma}_{n+1}^{(\ell)} =  \sum_l \widehat{\rho}_{n+1}^{(\ell)} S\widehat{\Sigma}_{f, n+1}^{(\ell)},$$
where $S = M S^n N^T$.
\item{\textbf{Truncation}} :  Determine the SVD $\widehat{S}^{n+1} = \widehat{P} \widehat{\Sigma} \widehat{Q}^T$ where $\widetilde{\Sigma} =  \text{diag}(\sigma_j)$. Choose 
$$r_{n+1} = r_n + \kappa \leq 2 r_n $$ 
Consider: $\widetilde{S}^{n+1} $ the diagonal matrix, whose diagonal is given by the vector of the $r_{n+1}$ largest singular values, $P_1 \in \mathbb{R}^{2r \times r_{n+1}}$ and $Q_1 \in \mathbb{R}^{2r \times r_{n+1}}$, the matrices made up by the first $r_{n+1} $ columns of $\widehat{P}$ and $\widehat{Q}$, respectively.\\
Finally, set $X^{n+1} = \widehat{U}P_1 \in \mathbb{R}^{m \times r_{n+1}}$ and $V^{n+1} = \widehat{V}Q_1 \in \mathbb{R}^{m \times r_{n+1}}$, 
\item Set $k_{n+1} = \| \widetilde{S}^{n+1} \| $ and $S^{n+1}= \widetilde{S}^{n+1} / k_{n+1}$, construct the approximation 
$$\phi ^{n+1} = X^{n+1} S^{n+1} V^{n+1, T}$$
which is the initial condition of rank $r_{n+1}$ to start DLRP.

\end{enumerate}

The previous iteration is inspired by the rank adaptive unconventional integrator of \cite{CKL},  adapted to the problem of \cite{KWMF}, using $(\phi_r, \keff^r)_{r=r_n}$ as initial condition and truncating at rank $r_{n+1} = r_n + \kappa$.

We consider two possibilities to choose the rank increase, i.e. considering $\kappa$ fixed or measuring the error with respect to the singular values $\widehat{\Sigma} =  \text{diag}(\sigma_j)$ as 
\begin{equation}
\label{rt}
 \sum_{j = r_{n+1}+1}^{2r} \sigma^2_j   \leq \varepsilon
\end{equation}
where $\varepsilon = \beta \Delta_n $, with $\Delta_n = | k_{n} - k_{n-1}|$ and $\beta<1$ a parameter  to measure how we truncate the singular values. 
The latter approach chooses the rank such that the error in the low-rank approximation is related to how close the power iteration is to convergence.

\begin{remark}
The introduction of proper rank control techniques for this type of problem is particularly important because our numerical simulations show that using the rank adaptive integrator constructed following \cite{CKL} and extended to this problem does not result in a convergent scheme. 
\end{remark}

\begin{remark}
The solution of the problem using the power method involves solving a linear matrix equation of size $N_x \cdot G \times N_x \cdot G$, of the form \eqref{eq_powerit}, at each step. This limits the possibility of solving such a system by direct methods if $N_x$ and $G$ are large. In contrast, the low rank power method in the $K, L$ and $S$ steps, involves solving matrix equations of size only $N_x \cdot r \times N_x \cdot r$, $G \cdot r  \times G \cdot r$ and $r^2 \times r^2$, respectively.  Since usually $r \ll \min \lbrace N_x, G \rbrace$, the low rank power method offers the possibility of storing the matrices involved and solving the systems by direct methods as well.

\end{remark}

\section{Numerical Experiments}

In this section, we present several numerical examples showing the effectiveness of the proposed rank adaptation algorithm. Note that full rank problems are considered in all experiments.
\subsection{Plutonium sphere}
\begin{figure}[h!]
\includegraphics[scale=0.4]{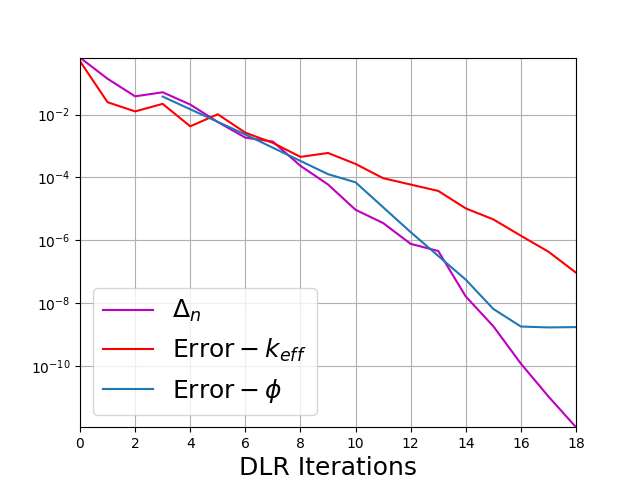}
\includegraphics[scale=0.4]{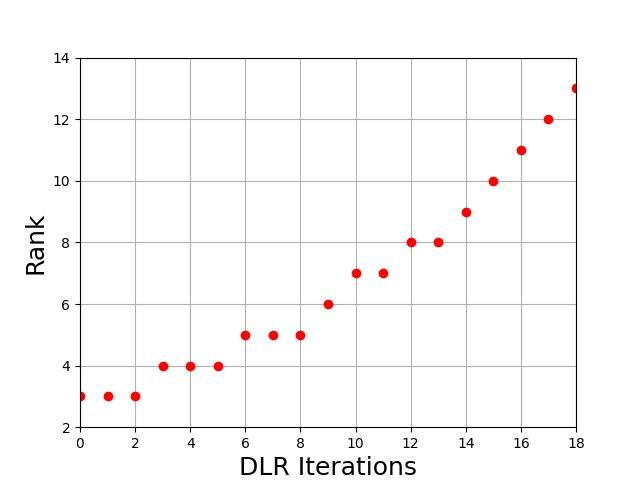}

\caption{Numerical results for a plutonium sphere, with fixed increase of rank $\kappa=1$. The left figure shows the behaviors of $\Delta_n$ and of the errors in both the eigenvalue and the eigenmatrix. In the right plot, the adaptation of the rank is presented. }\label{plutonium}
\end{figure}

As first example, we consider a plutonium sphere.
The energy domain is represented by $70$ energy groups and a  spatial discretization with $70$
spatial cells is chosen. 
The reference solution has $\keff^*=  0.9956975948475687$, obtained by solving the full rank problem by the power method.
The starting tolerance is set to $\theta_0 = 0.1$ and the final tolerance is $\theta = 10^{-12}$.
We set $\rho = 0.1$ and a fixed increment for the rank, $\kappa = 1$.
The initial rank is $r_0= 5$ and we finish with rank $r= 13$, after $19$ iterations with  $\text{error} = 4.25 \times 10^{-9}$ (w.r.t. to the eigenvalue).
The results of the simulation are shown in Figure \ref{plutonium}.
In this problem a rank of $13$ is required to satisfy the desired tolerance. Without rank adaption the entire iterations would need to be run using this rank. Here, however, the iterations are conducted using an average rank of only $6.8$. This results in a significant reduction in computational cost (note that computational cost scales quadratic in rank).

\subsection{Stainless-steel reflected uranium sphere}
A standard benchmark problem is the IEU-MET-FAST-005 criticality benchmark from the OECD/NEA suite, see \cite{Gorbatenko}.
This problem has a sphere of $36\%$ enriched uranium surrounded by a neutron reflector comprised of stainless steel. The problem has
an overall radius of $21.486$ cm and the uranium sphere has a radius of $13.213$ cm. The stainless steel is divided into shells
with two different densities: one with radius $1.758$ cm and the other covers the remainder of the total size.
 A spatial discretization with $400$
spatial cells is considered and the energy domain is represented by $87$ energy groups. 
The reference solution can be computed by the power method  and it is given by $\keff^*=0.9269867446024421$.

\begin{figure}[h!]
\includegraphics[scale=0.4]{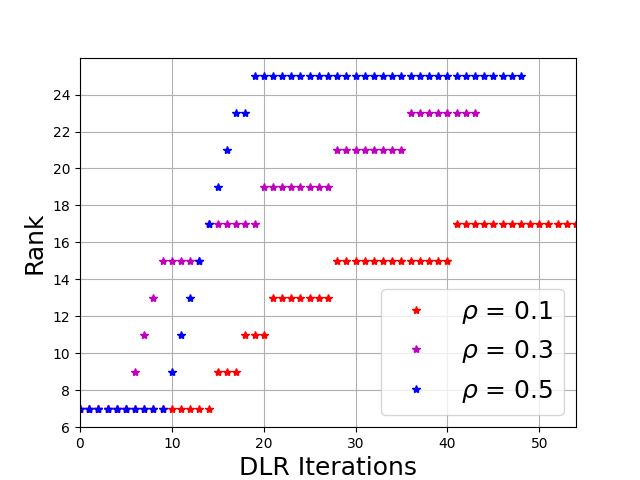}
\includegraphics[scale=0.4]{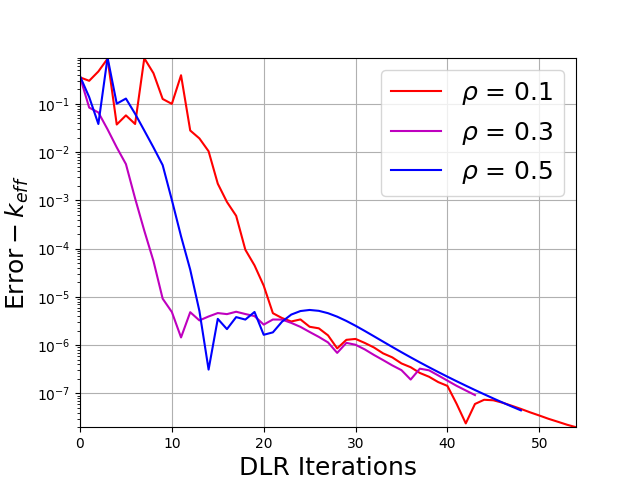}
\includegraphics[scale=0.4 ]{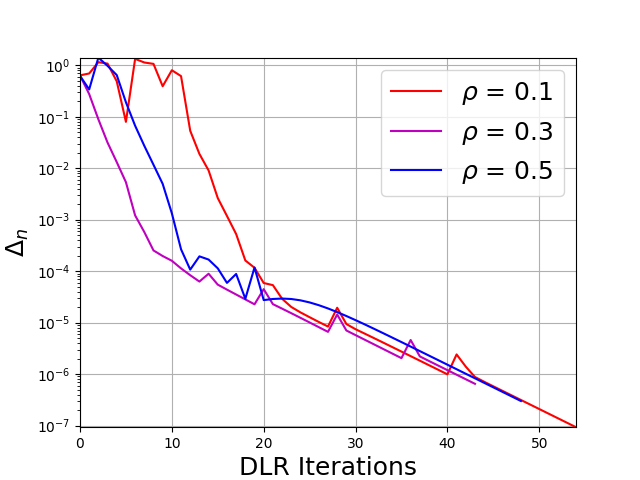}
\includegraphics[scale=0.4]{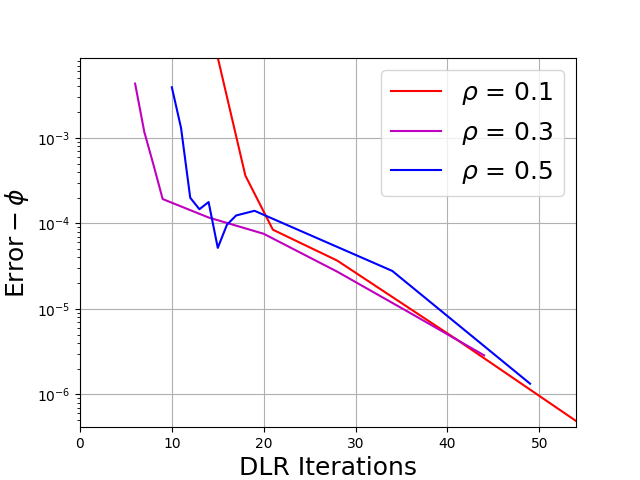}
\includegraphics[scale=0.4]{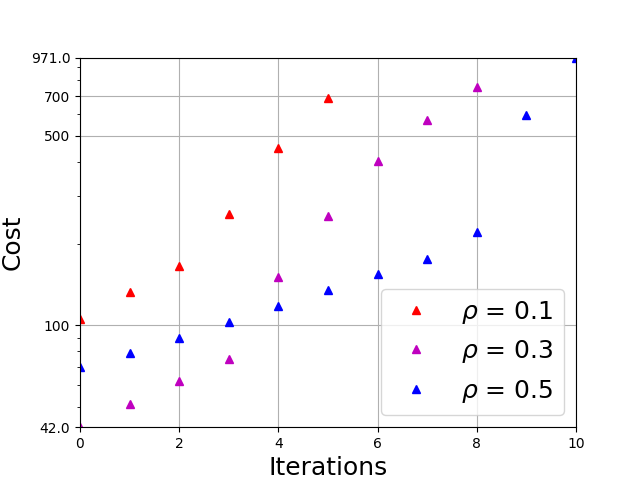}
\caption{Numerical results of stainless-steel reflected uranium sphere with the rank adaptive algorithm for different values of $\rho$.}\label{stainless_all}
\end{figure}

We consider $\theta_0 = 10^{-2}$, $\theta= 10^{-7}$ and $\kappa = 2$.
The initial rank is $r_0=7$.
The behavior of the simulations is shown in Figures \ref{stainless_all}, for three different values of $\rho$, the factor of reduction of the tolerance.
In particular, we show the errors and the rank adaptation with respect to the DLR iterations. 
The last plot in Figure \ref{stainless_all} shows the computational cost with respect to the iterations of Algorithm \ref{algorithm1}.
The cost of any iteration is computed as the sum of all the DLR iterations, weighted by the corresponding rank.
The results show that for a wide range of choices for $\rho$ the algorithm works well and gives accurate results. However, for larger $\rho$ the rank increases relatively quickly, which requires the algorithm to perform more iterations with higher rank (negatively impacting computational cost). Overall this results in a performance difference of approximately 40\% between the values of $\rho$ studied here.

\begin{figure}[h]
\includegraphics[scale=0.4]{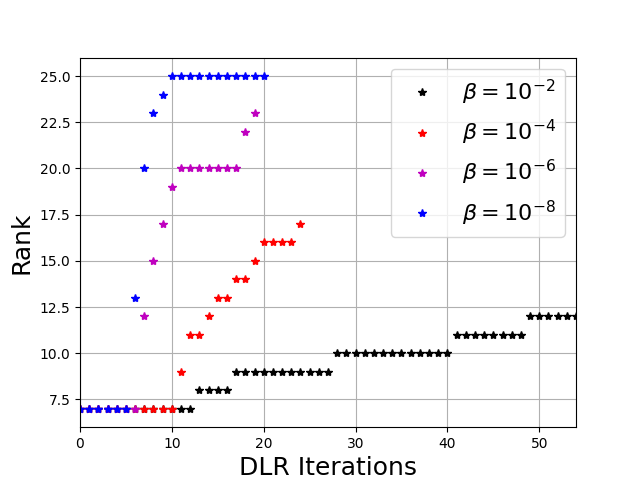}
\includegraphics[scale=0.4]{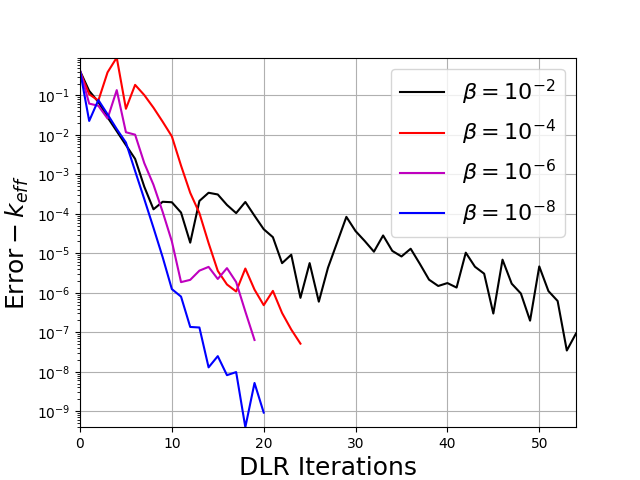}
\includegraphics[scale=0.4]{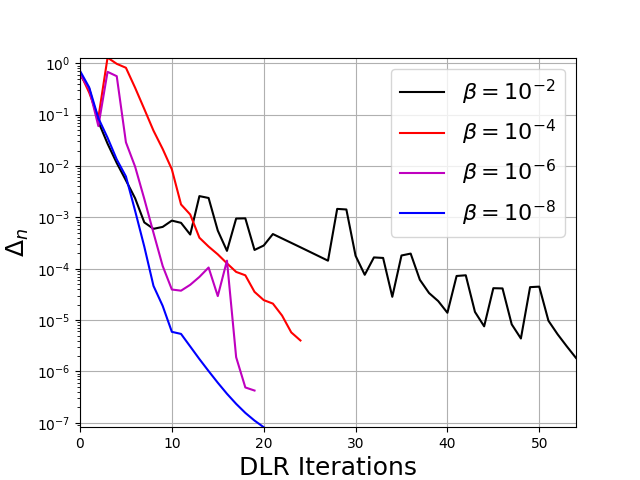}
\includegraphics[scale=0.4]{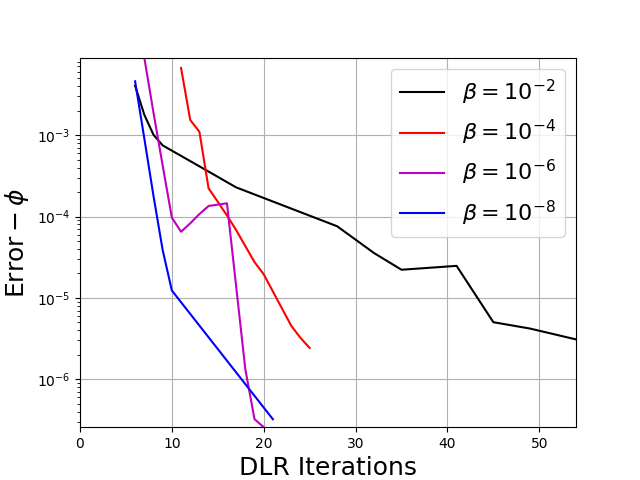}
\includegraphics[scale=0.4]{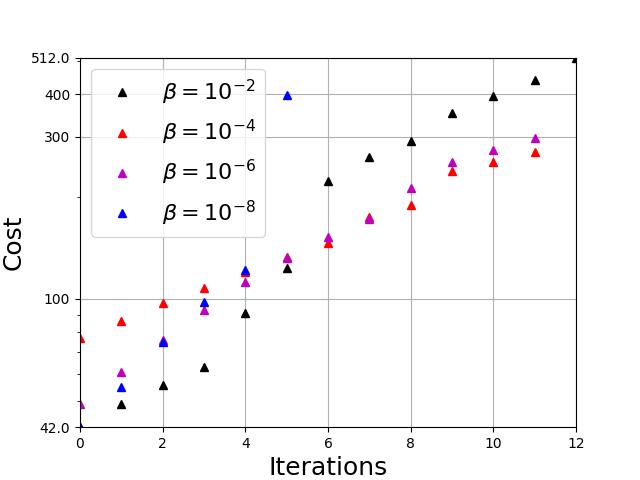}
\caption{Numerical results of stainless-steel reflected uranium sphere, with adaptation with some values of $\beta$ and $\rho= 0.5$.}\label{stainless_all1}
\end{figure}

We also observe (see Figure \ref{stainless_all1}) that the algorithm works well over a wide range of $\beta$. For $\beta = 10^{-8}$, which allows a fast increase of the rank, we get very accurate results, but also a significantly higher computational cost in comparison to $\beta = 10^{-4}$ and $\beta = 10^{-6}$. For larger values of $\beta$, here we consider $\beta = 10^{-2}$, the rank increases slowly, also leading to a higher  computational cost. But for the values in between (we choose $\beta=10^{-4}$ and $\beta = 10^{-6}$) we get accurate results at low computational cost.
We note that we enforce here that the rank is monotonically increasing. This is only necessary for large $\beta$ where in some situations the algorithm would otherwise decrease the rank (which is clearly not desirable).  

\subsection{Light water reactor}

In this problem we look at the solution for a problem of a homogenized light water reactor using the SHEM
$361$ energy group structure, \cite{hebert2008}.
 The problem consists of this homogenized material in a sphere of
radius $79.06925$ cm. 
The reference solution is given by $\keff^*= 0.9999061310852358 $.  A spatial discretization with $400$
spatial cells is chosen. The energy domain is represented by $361$ energy groups. 

\begin{figure}[h]
\includegraphics[scale=0.4]{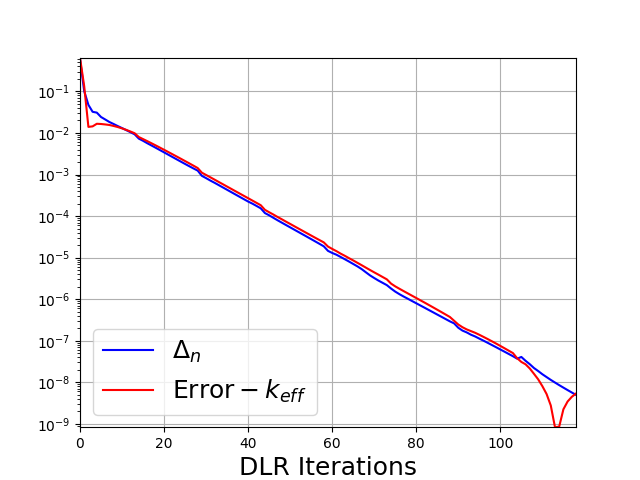}
\includegraphics[scale=0.4]{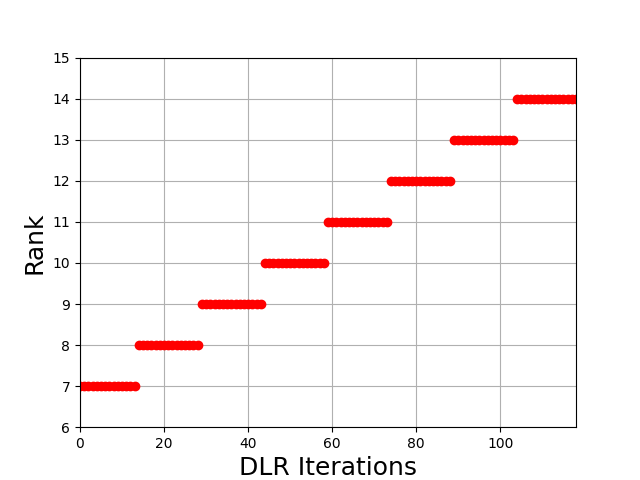}
\caption{Numerical results of light water ractor.}\label{water}

\end{figure}

In the numerical simulation, the starting rank is $r_0=7$ and  we use a fixed  increment of the rank $ \kappa = 1$.
The tolerances are $\theta= 10^{-9}$ and $\theta_0 = 10^{-2}$. The parameter $\rho$ is set to $0.25$.
The results  are showed in Figure \ref{water}.

\section{Optimization}
Computing the effective eigenvalue for a fixed configuration, as has been done up this point, is useful in physical applications. However, often we are interested in modifying a configuration (specified, e.g., by a parametrization of the geometry of the reactor) such that a certain objective is achieved.

In particular, here we are interested in computing the values for the parameters such that the resulting system is critical. We consider the following inverse problem: for a given target value of the effective eigenvalue $\keff^*$ and a  physical system, parameterized by $\alpha_1, ..., \alpha_n$, determine values of the parameters such that, the corresponding effective eigenvalue coincides with $\keff^*$.
We can formulate an optimization problem, i.e. compute 
\begin{equation}
\label{optpb}
 \;\;\alpha_1^*, ..,\alpha_n^* = \argmin_{\alpha_1, ..,\alpha_n \geq 0} f(\alpha_1, ..,\alpha_n)
\end{equation}
where $f(\alpha_1, ..,\alpha_n) = | \keff(\alpha_1, ..,\alpha_n)-\keff^*|$.\\
When $\keff^*=1$, we are optimizing the distance to criticality.
The idea is to solve \eqref{optpb} by a low-rank based gradient algorithm.
This assumes the computation of $f(\alpha_1, ..,\alpha_n)$ and, consequently, of $\keff(\alpha_1, ..,\alpha_n)$, at any step of the algorithm (and in the approximation of the gradient). Computational efficiency in this context is even more important as we operate in a multi-query context (i.e.~the optimization algorithm needs to solve the eigenvalue problem many times). The possibility of computing the objective function using the dynamical low rank power method can thus be crucial.

For a given parameterized system and fixing the rank $r$ a priori, the resulting optimization procedure is presented in algorithm \ref{algorithm2}.  
We refer to \cite{NW} for the gradient algorithm and its variants.
We have that $\alpha_1^0$,.., $\alpha_n^0$ are the starting values for the parameters of the physical model and $ \phi_0 = U_0 S_0 V_0^T$ is the starting randomly chosen eigenmatrix of rank $r$.  The objective function is computed by the low rank power method with accuracy $tol_f$.
The parameter $tol$ is the total accuracy required for the optimization.

 \begin{algorithm}[h!]
 \caption{\textbf{Low rank optimization}}\label{algorithm2}
 \KwData{ $\alpha_1^0$,.., $\alpha_n^0$, $ \phi_0 = U_0 S_0 V_0^T$, $tol$, $tol_f$, $h$, $c$, $p$, $r$}
    
  $f_1 = f(\alpha_1^0$,.., $\alpha_n^0, \phi_0 = U_0 S_0 V_0^T, r, tol_f)$\
    
    \While {$f_1> tol$ $\&$ $h>h_{min}$}{
     $\nabla f_0 = \nabla f(\alpha_1^0$,.., $\alpha_n^0, \phi_0,  tol_f)$  \Comment*[r]{ via finite differences}\
     $ \alpha_i = \alpha_i^0 - h \nabla_{\alpha_i} f_0, \quad \forall i = 1, ...,n$ \Comment*[r]{Gradient step}\
     
         $\left[ f_1, \phi=  USV^T \right]  = f(\alpha_1$,.., $\alpha_n, \phi = U_0 S_0 V_0^T,  tol_f)$\\

        \While{$f_1 > f_0-c h \| \nabla f_0 \|^2$  \Comment*[r]{Armijo condition}}{ 
         $h = ph$ \Comment*[r]{Stepsize reduction}
        $ \alpha_i = \alpha_i^0 - h \nabla_{\alpha_i} f_0, \quad \forall i = 1, ...,n$ \Comment*[r]{Gradient step}
     
       $\left[ f_1, \left( \phi=  USV^T\right)  \right]  = f(\alpha_1,.., \alpha_n, \phi_0 = U_0 S_0 V_0^T,  tol_f)$}
       
    $f_0= f_1$;\\
     $\alpha_i^0 = \alpha_i, \quad \forall i = 1,...,n$;\\
    $\phi_0 = \phi$}
   
\end{algorithm}

A further possibility is to include the previous low rank optimization in an adaptive rank setting. The adaptive scheme is presented in the algorithm \ref{algorithm3}.
We start with a (small) rank $r_0$ and proceed with an optimization step, performed by the previous algorithm \ref{algorithm2}, up to a (low) desired accuracy. 
Then, we reduce the tolerance by a factor $\rho<1$ and we use the computed values of the parameters to start a new optimization step  with a higher rank $r$ to get a more accurate solution.
 In this way, we obtain a more accurate solution  with respect to the rank truncation error of the eigenmatrix. Clearly, to start the optimization with higher rank $r$ we need to merge the previuos solution of rank $r_0$, in the manifold $\mathcal{M}_r$ via the function CR introduced in the previous section.  The advantage here is that when we are far away from criticality only a small rank is used. This is computationally efficient. As we get closer to criticality, we increase the rank to get the desired accuracy.

 \begin{algorithm}[h!]
 \caption{\textbf{Rank adaptive optimization}}\label{algorithm3}
     \KwData{$\alpha_1^0,.., \alpha_n^0$, $\phi_0 = U_0 S_0 V_0^T $, $tol$, $tol_f$, $r_0$, $\kappa$, $\rho$}
 $r = r_0$ \
    $\left( \alpha_1,.., \alpha_n^, f,\phi_r = U_r S_r V_r^T \right)\longleftarrow \text{optimize}(\alpha_1^0$,.., $\alpha_n^0, \phi_0, tol_f, r) $ \
   \While{$f <tol$} {
   $\phi_{r+ \kappa} = \text{CR}(\phi_r)$ \Comment*[r]{merge the  solution in $\mathcal{M}_{r+\kappa}$}\
 $r = r+\kappa$  \Comment*[r]{increase the rank}\
    $\alpha_i = \alpha_i^0 \quad \forall i = 1, ...,n$ \Comment*[r]{Re-initialization}\
   $\phi_0 = \phi_{r+ \kappa}$\\
  $tol =  \rho tol$\
  $\left( \alpha_1,.., \alpha_n,f, \phi_r = U_r S_r V_r^T\right)\longleftarrow \text{optimize}\left( \alpha_1^0,.., \alpha_n^0, \phi_0, tol_f, r \right)  $ }

\end{algorithm}

\section{Numerical experiments}

In this section, we test the low rank optimization methods presented in the previous section on a physical test problem.
We have a sphere with a hollow center, surrounded by uranium U, and then two types of stainless steel, denoted by SS2 and SS3 (i.e. 4 layers in total). 
In principle, the length of any of those layers (i.e. length hollow, lengthU, lengthSS2, lengthSS3) could be a parameter for optimization. 
In this case, we consider lengthU = $\alpha_1$ and lengthSS2 = $\alpha_2$ as parameters. The only constraint on $\alpha_1$ and $\alpha_2$ is positivity. 

In Table \ref{grad_test_fr}, we show the results obtained for the fixed rank optimization algorithm applied to the before mentioned problem, for different values of the rank. 
In the experiments, we choose $c= 10^{-4}$, $p= 0.5$, $h_{min} = 10^{-9}$, $tol_f = 10^{-7}$ and $tol = 10^{-7}$. 
The starting values are $\alpha_1^0=  1.75$ and $\alpha_2^0 = 10.564$. 
We show the computed optimal values for the parameters $\alpha_1$ and $\alpha_2$, the objective function, the total computational cost and the low-rank rank error.
 In particular, the computational cost for the gradient iterations in Table \ref{grad_test_fr} is  the sum of all the low rank power method iterations done. The cost also includes the approximation of the gradient and of the line search.
 For this problem, for any $\alpha_1$ and $\alpha_2$ it is possible to compute the effective eigenpair, and thus the objective function using the full power method.
Therefore, for any gradient iteration we are able to evaluate the rank truncation error of the eigenmatrix, thanks to the comparison of the one computed by the full rank solver.
The rank error presented in the Table \ref{grad_test_fr}, is the  error between the full and the low-rank eigenmatrices.

\begin{figure}[h]

 \centering
    \includegraphics[scale=0.3]{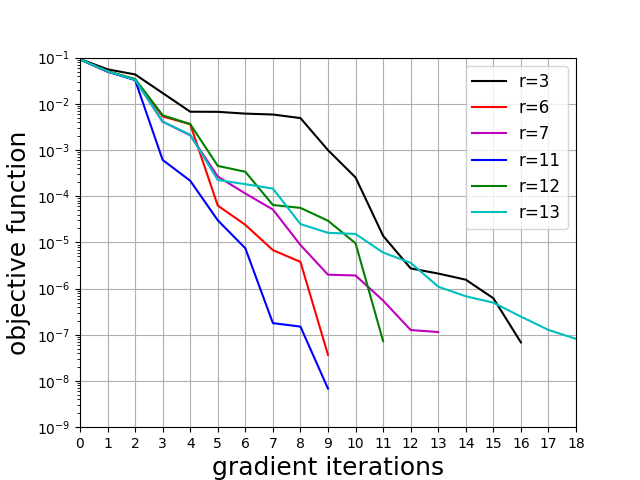}
     \includegraphics[scale=0.3]{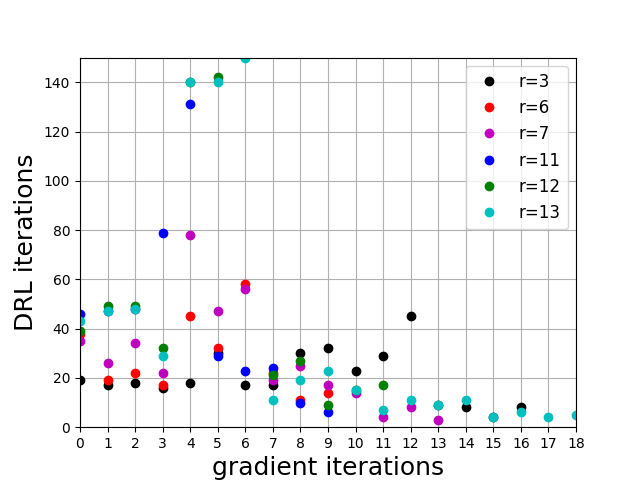}
     \includegraphics[scale=0.3]{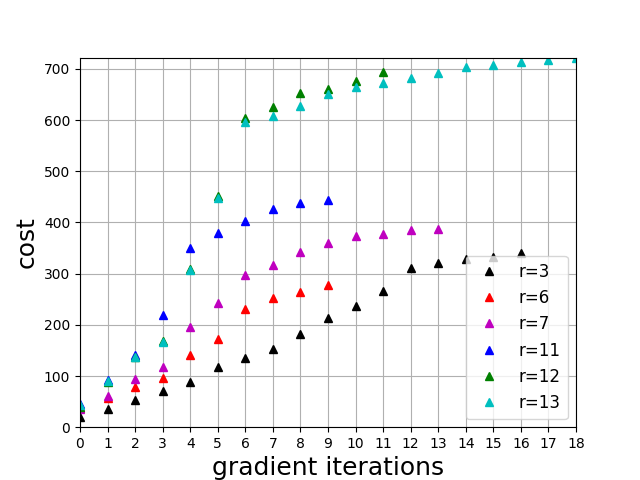}
   
    \caption{Computational results for the fixed rank gradient optimization.  The computational cost of any gradient step is computed as the sum of all the previous low-rank power method iterations.}
    \label{grad_plot_fr}
\end{figure}

\begin{figure}[h]
\centering
     \includegraphics[scale=0.3]{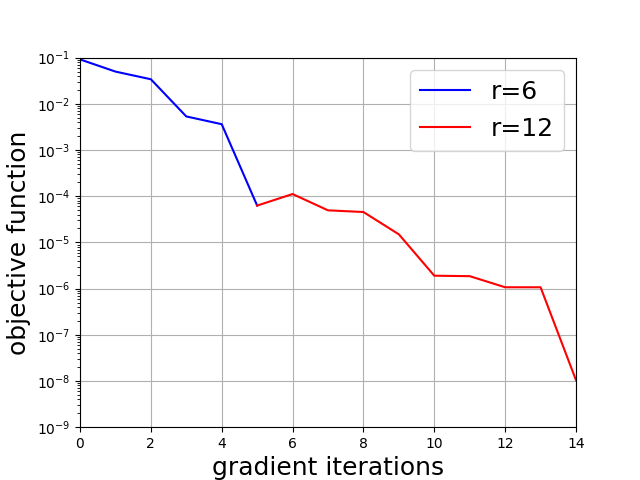}
     \includegraphics[scale=0.3]{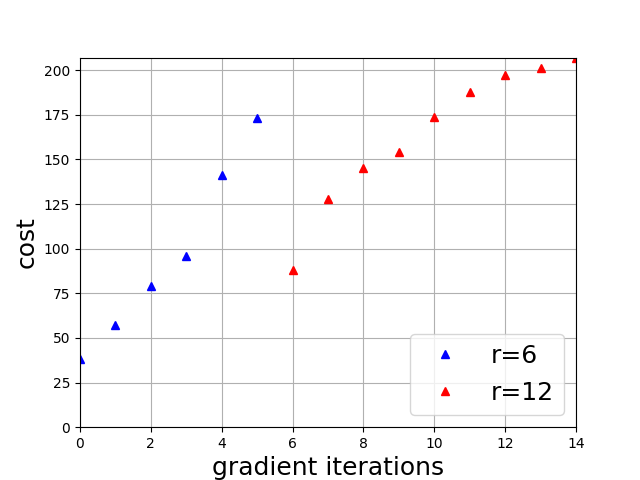}
    \includegraphics[scale=0.3]{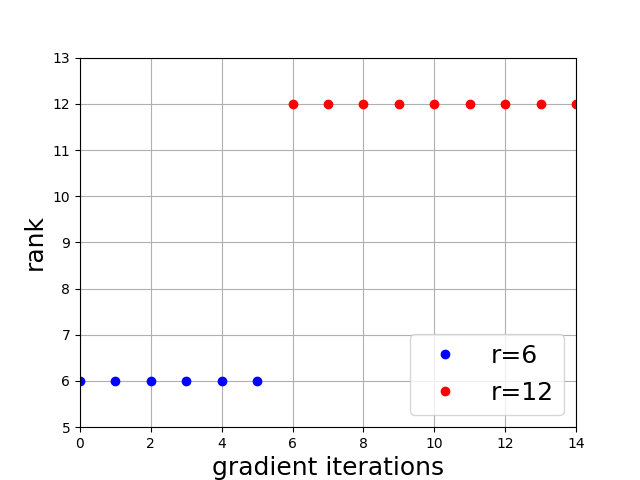}
    \includegraphics[scale=0.3]{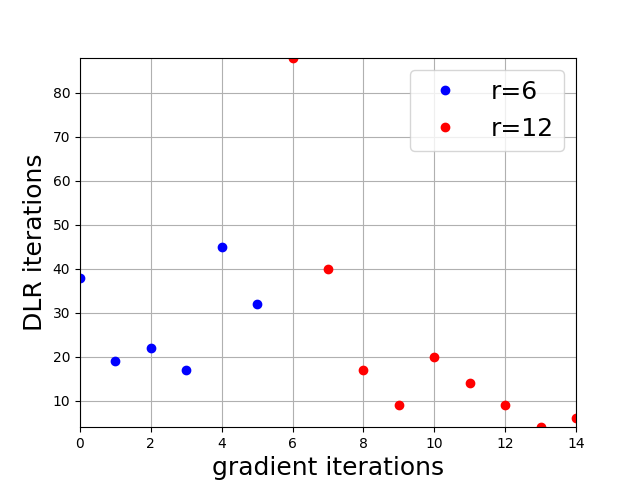}
    \caption{Computational results for rank adaptation. The starting rank is $6$ and then it is increased to $12$. The computational cost of any gradient step is computed as the sum of all the previous low-rank power method iterations.}
    \label{612}
\end{figure}

\begin{figure}[h]
\centering
    \includegraphics[scale=0.3]{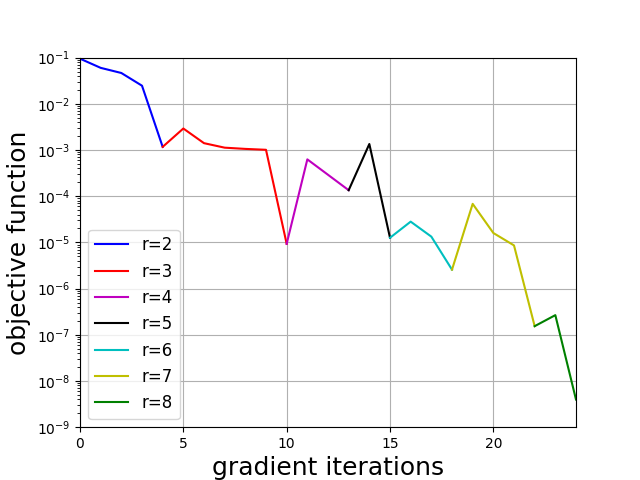}
     \includegraphics[scale=0.3]{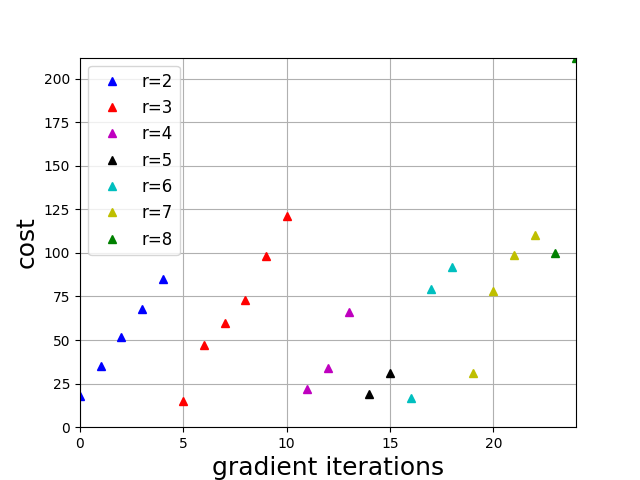}
     \includegraphics[scale=0.3]{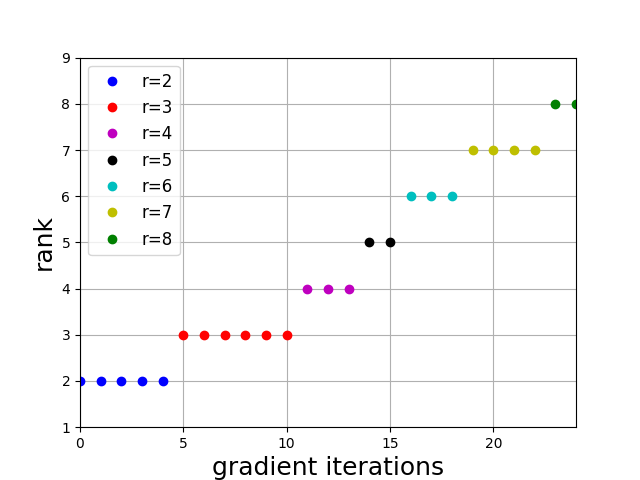}
     \includegraphics[scale=0.3]{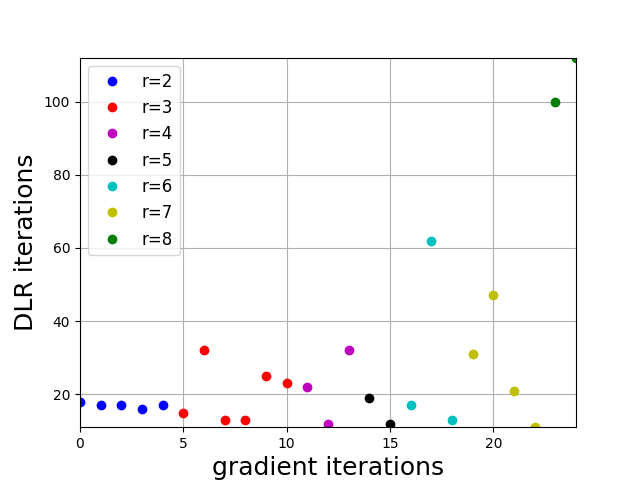}
    \caption{Computational results for rank adaptation. The starting rank is $2$. The computational cost of any gradient step is computed as the sum of all the previous low-rank power method iterations.}
    \label{2to8}
\end{figure}

The computational results in Table \ref{grad_test_fr} and Figure \ref{grad_plot_fr} show that, for any fixed value of the rank, we are able to obtain an approximation of the eigenvalue with good accuracy. However, it is necessary to consider  moderately higher ranks to obtain good accuracy, with respect to the eigenmatrix, in terms of the rank truncation error.
This motivates the use of a rank adaptive optimization.
In Figures {\ref{612} and \ref{2to8}}, we show the results of the adaptive optimization.
In particular, in the experiment of Figure {\ref{612}}, we start with a rank $6$ and a first stopping tolerance of $10^{-4}$. Then, we reduce the tolerance to $10^{-7}$ and directly double the rank. 
At the end, we get the same order of accuracy for the objective function with respect to the simulation with fixed rank $12$. 
However, we have a total number of iterations of $376 $, of which $ 173 $ are done with rank $6$.
Thus, we have a considerable saving in computational cost if we compare to the $ 693$ iterations at fixed rank $12$.

In the test of {Figure \ref{2to8}}, we show that it is also possible to gradually increase of the rank, starting with rank $2$ and tolerance $10^{-2}$, reducing at each step by a factor $\rho = 0.15$.  We finally obtain a rank $8$ eigenmatrix.

The final error for $\keff$ is $3.97 \times 10^{-9}$. The number of iterations with the largest rank has been significantly decreased.

To show the importance  of a sufficient rank increase, we solve by the power method the problems obtained considering the optimal values of the parameters computed in the previous adaptive experiments. 
In particular, when we solve with the parameters obtained from the experiment with final rank $12$, we obtain an eigenvalue that differs from $1$ by an error of $2.349277 \times 10^{-8}$.
The same computation with the values of $\alpha_1$ and $\alpha_2$ from the second experiment with final rank $8$, gives an error of $2.303825 \times 10^{-5}$.

\begin{table}[h]                                                                         
\centering                                                                            
\begin{tabular}{llllll}                                                            
\hline
rank & $f(\alpha_1, \alpha_2)$   & $\alpha_1$ & $\alpha_2$ & cost & rank error\\
\hline
$3$ &  $  6.852542 \times 10^{-8}  $  &  $ 12.438101 $  &   $ 3.550126  $  &   $ 340  $ &    $1.840924 \times 10^{-3}    $  \\
\hline
  $6$ &  $ 3.687297 \times 10^{-8}   $  &  $ 12.360320 $  &   $  3.606547  $  &   $ 278   $   &   $ 4.858112 \times 10^{-5}  $\\
  \hline
  $7$ &  $ 1.147334 \times 10^{-7}   $  &  $ 12.361202 $  &   $ 3.554740 $  &   $ 388   $  &   $ 2.225781 \times 10^{-5}   $\\
  \hline
  $11$ &  $ 2.458198\times 10^{-8}   $  &  $ 12.359178 $  &   $ 3.555892  $  &   $ 443 $\  &  $ 3.755094 \times 10^{-7}  $\\
  \hline
  $12$ & $7.277437 \times 10^{-8}$ &  $  12.308510
 $  & $  3.606547 $  & $693$  &  $  1.963992 \times 10^{-7}  $\\
  \hline
  $13$ & $   7.906809 \times 10^{-8}$ & $ 12.347656$  & $ 3.567397  $ &  $ 729 $   & $ 3.582495 \times 10^{-8}  $\\                                                                    
\hline
\end{tabular}
\caption{Computational results for the fixed rank gradient optimization.}
\label{grad_test_fr}
\end{table}

\begin{table} 
\label{ad1}                                                                        
\centering                                                                            
\begin{tabular}{llllll}                                                            
\hline
 $\alpha_1$ & $ \alpha_2$ & objective function & error & $ r_0$ & $ r_f$\\
 
 \hline
$ 12.361255$  & $3.553801$ &  $ 1.150259\times 10^{-7}$  & $2.348753\times 10^{-8}   $ & $ 6$ & $12$\\
\hline
$ 12.421492 $  &  $ 3.494477 $    & $ 3.968996 \times 10^{-9} $ &  $2.303825 \times 10^{-5}  $ & $2$ & $8  $\\
\hline
\end{tabular}
\caption{Final results for the experiments with rank adaptation. $r_0$ and $r_f$ are, respectively, the initial and the final ranks. The error is the the distance to $1$ of the eigenvalue computed by the power method.}
\end{table}

\begin{acknowledgements}
The authors thank Alexander Ostermann (University of Innsbruck) for helpful discussions.
C. Scalone  thanks the INDAM
Research group GNCS for financial support. 
\end{acknowledgements}

%
 \section*{Conflict of interest}
 The authors declare that they have no conflict of interest.


\newpage
\begin{thebibliography}{}
\bibitem{Buchan}A. G. Buchan, C. C. Pain, F. Fang and  I. M. Navon, A POD reduced-order model for eigenvalue problems with application to reactor physics, Int. J. Num. Methods Engineering Volume95, 12 (21), 1011-1032 (2013).
\bibitem{lrkin4} F. Cassini, L. Einkemmer, Efficient 6D Vlasov simulation using the dynamical low-rank framework Ensign,  Comput. Phys. Commun.,  280:108489, (2022). 
\bibitem{CL} G. Ceruti, C. Lubich, An unconventional robust integrator for dynamical low-rank approximation, BIT Numer. Math., 62  23-44 (2022).
\bibitem{CKL} G. Ceruti, J. Kusch, C. Lubich, A rank-adaptive robust integrator for dynamical low-rank approximation, BIT Numer. Math., 62, 1149–1174 (2022).
\bibitem{lradap3} G. Ceruti, J.  Kusch, C. Lubich, A parallel rank-adaptive integrator for dynamical low-rank approximation, SIAM J. Sci. Comput., 46(3), B205--28 (2024). 
\bibitem{lrap4} J. Coughlin, J.  Hu, Efficient dynamical low-rank approximation for the Vlasov-Ampere-Fokker-Planck system. J. Comput. Phys., 470, 111590, (2022). 
\bibitem{lrcons4} J. Coughlin, J.  Hu , U. Shumlak, Robust and conservative dynamical low-rank methods for the Vlasov equation via a novel macro-micro decomposition, J. Comput. Phys., 509, 113055,  (2024). 

\bibitem{lrap1} Z. Ding, L. Einkemmer, Q.  Li, Dynamical low-rank integrator for the linear Boltzmann equation: error analysis in the diffusion limit. SIAM J. Numer. Anal., 59(4), 2254--2285, (2021). 
\bibitem{lrap3}  L. Einkemmer, J. Hu, J. Kusch, Asymptotic-Preserving and Energy Stable Dynamical Low-Rank Approximation, SIAM J.  Numer.Anal., 62(1), 73--92, (2024). 
\bibitem{lrap2}L. Einkemmer L, J. Hu, Y. Wang, An asymptotic-preserving dynamical low-rank method for the multi-scale multi-dimensional linear transport equation,  J. Comput. Phys., 439, 110353, (2021).

\bibitem{lrkin1} L. Einkemmer, C. Lubich, A low-rank projector-splitting integrator for the Vlasov--Poisson equation, SIAM J. Sci. Comput., 40(5), 1330--1360 (2018).
\bibitem{lrbiology3} L. Einkemmer, J. Mangott, M. Prugger, A low-rank complexity reduction algorithm for the high-dimensional kinetic chemical master equation. J. Comput. Phys., 6, 112827, (2024).

\bibitem{lrcons1} L. Einkemmer, I. Joseph, A mass, momentum, and energy conservative dynamical low-rank scheme for the Vlasov equation, J. Comput. Phys.,  443,110495, (2021).
\bibitem{lrkin2} L. Einkemmer, A. Ostermann, C. Piazzola, A low-rank projector-splitting integrator for the Vlasov–Maxwell equations with divergence correction,  J. Comput. Phys., 403, 109063, (2020).

\bibitem{EOS}L. Einkemmer, A. Ostermann, C. Scalone, A robust and conservative dynamical low-rank algorithm, Journal of Computational Physics, (2023).

\bibitem{GR}P. German, J.C. Ragusa
Reduced-order modeling of parameterized multi-group diffusion k-eigenvalue problems,
Ann. Nucl. Energy, 134, 144-157 (2019).
\bibitem{GKS} N. Guglielmi, D. Kressner, and C. Scalone, Computing low-rank rightmost eigenpairs of a class of matrix-valued linear operators, Adv. Comp. Math., 47(5), 1–28, 2021.
\bibitem{Gorbatenko}M. V. Gorbatenko, V. P. Gorelov, V. P. Yegorov, V. G. Zagrafov, A. N. Zakharov, V. I. Ilyin, M. I. Kuvshinov, A. A. Malinkin, and V. I. Yuferev. Steel-reflected spherical assembly of 235U(36 $\%$). In International Handbook of Evaluated Criticality Safety Benchmark Experiments, number NEA/NSC/DOC/(95)03/III. Organisation for Economic Co-operation and Development Nuclear Energy Agency, September 2002.
\bibitem{HS} C. D. Hauck and S. Schnake, A Predictor-Corrector Strategy for Adaptivity in Dynamical Low-Rank Approximations, SIAM J. Matrix Anal. Appl., 44(3), 971-1005 (2023).
\bibitem{hebert2008} A. Hébert, A. Santamarina, Refinement of the Santamarina-Hfaiedh energy mesh between 22.5 eV and 11.4 keV. In International Conference on the Physics of Reactors, Interlaken, Switzerland, (2008).
\bibitem{HNS} M. Hochbruck, M. Neher and S.Schrammer, Rank-adaptive dynamical low-rank integrators for first-order and second-order matrix differential equations, BIT Numer. Math., 63, 9, (2023). 
\bibitem{lrbiology1} T. Jahnke, W. Huisinga, A dynamical low-rank approach to the chemical master equation, Bull. Math. Biol., 70, 2283--22302, (2008).
\bibitem{KPN} D. Knoll, H. Park, C. Newman, Acceleration of k-eigenvalue/Criticality calculations csing the Jacobian-free Newton-Krylov method, Nuclear Science and Engineering, 167, 133-140 (2017). 
\bibitem{lrkin5} J. Kusch, P. Stammer, A robust collision source method for rank adaptive dynamical low-rank approximation in radiation therapy,  ESAIM: Math. Model. Numer. Anal., 57(2), 865--891 (2023). 
\bibitem{uncertainty1} J. Kusch, G. Ceruti, L. Einkemmer, M. Frank, Dynamical low-rank approximation for Burger's equation with uncertainty, Int. J. Uncertain. Quantif., 12(5), (2022).
\bibitem{KL}O. Koch and C. Lubich, Dynamical low-rank approximation. SIAM J. Matrix Anal. Appl., 29(2), 434–454 (2007).
\bibitem{KWMF} J. Kusch, B. Whewell, R. McClarren and M. Frank, A low-rank power iteration scheme for neutron transport criticality problems, J. Comput. Phys., 470,  111587 (2022).
\bibitem{lrqm2} C. Lubich, From quantum to classical molecular dynamics: reduced models and numerical analysis. European Mathematical Society, (2008).
\bibitem{LO} C. Lubich and I. V. Oseledets, A projector-splitting integrator for dynamical low-rank approximation, BIT Numer. Math., 54(1), 171–188, (2014).


\bibitem{McClarrenBook}R. McClarren. Computational nuclear engineering and radiological science using python. Academic Press, 2017.
\bibitem{M} R. G. McClarren, Theoretical aspects of the simplified Pn equations, Transport Theor. Stat., 39(2-4), 73–109 (2010).
\bibitem{lrqm1} H.D. Meyer, F.Gatti, G. A. Worth, editors. Multidimensional quantum dynamics: MCTDH theory and applications. John Wiley \& Sons, (2009).

\bibitem{lrqm3} A. Nonnenmacher, C. Lubich, Dynamical low-rank approximation: applications and numerical experiments, Math. Comput. Simulat., 79(4), 1346--1357 (2008). 
\bibitem{NW} J. Nocedal, S. J. Wright, Numerical Optimization, 2nd Edition, Springer, 2006.
\bibitem{Ortega} M.I. Ortega , R.N. Slaybaugh , P.N. Brown , T.S. Bailey , B. Chang, A Rayleigh quotient method for criticality eigenvalue problems in neutron transport,  Ann. Nucl. Energy, 138, 107120 (2020).
\bibitem{PMF} Z. Peng, R. McClarren, and M. Frank, A low-rank method for two-dimensional time-dependent radiation transport calculations, J. Comput. Phys., 421(109735), (2020).
\bibitem{PM} Z. Peng and R. McClarren, A high-order/low-order (HOLO) algorithm for preserving conservation in time-dependent low-rank transport calculations, J. Comput. Phys., 447(110672)  (2021).
\bibitem{lrbiology2} M. Prugger, L. Einkemmer, C. F. Lopez, A dynamical low-rank approach to solve the chemical master equation for biological reaction networks,  J. Comput. Phys., 489, 112250, (2023).
\bibitem{StaceyBook} W. M. Stacey. Nuclear reactor physics. John Wiley $\&$ Sons, 2018.
\bibitem{lrqm4} D. Sulz, C. Lubich, G. Ceruti, I. Lesanovsky, F. Carollo, Numerical simulation of long-range open quantum many-body dynamics with tree tensor networks, Phys. Rev. A, 109(2) (2024). 
\bibitem{Sun}Y. Sun, Y. Yang , Y. Wang , Z. Li , Y. Ma,  A POD reduced-order model for resolving the neutron transport problems of nuclear reactor, Ann. Nucl. Energy, 149 (107799) (2020).
\bibitem{WPK} J. Willert, H. Park, D.A. Knoll, A comparison of acceleration methods for solving the neutron transport k-eigenvalue problem, J. Comput. Phys., 274, 681-694 (2014).



























\end{thebibliography}


\end{document}